\newtheorem{lemma}{Lemma}[section]
\newtheorem{prop}[lemma]{Proposition}
\newtheorem{theorem}[lemma]{Theorem}
\newtheorem{cor}[lemma]{Corollary}
\newtheorem{rem}[lemma]{Remark}
\newtheorem{exam}[lemma]{Example}

\newcommand{\kla}{\left ( }
\newcommand{\nach}{\rightarrow}
\newcommand{\mer}{\right ) }

\newcommand{\for}{\begin{eqnarray*}}
\newcommand{\mel}{\end{eqnarray*}}

\newcommand{\mitt}{\left | { \atop } \right.}
\newcommand{\kl}{\pl \le \pl}

\newcommand{\lel}{\pl = \pl}
\newcommand{\lell}{\p=\p}
\newcommand{\kll}{\p \le \p}
\newcommand{\gll}{\p \ge \p}

\newcommand{\ov}{\overline}

\newcommand{\nz}{{\rm  I\! N}}
\newcommand{\nen}{n \in \nz}
\newcommand{\ken}{k \in \nz}
\newcommand{\ez}{{\rm I\! E}}
\newcommand{\dz}{ {\rm I\!D}}

\newcommand{\cz}{{\rm I\!\!\! C}}

\newcommand{\p}{\hspace{.05cm}}
\newcommand{\pl}{\hspace{.1cm}}
\newcommand{\pll}{\hspace{.3cm}}
\newcommand{\pla}{\hspace{1.5cm}}
\newcommand{\hz}{\vspace{0.5cm}}

\newcommand{\al}{\alpha}

\newcommand{\si}{\sigma}

\newcommand{\la}{\lambda}
\newcommand{\eps}{\varepsilon}

\newcommand{\ds}{D_{\si}}

\newcommand{\lzn}{\ell_2^n}

\newcommand{\lin}{\ell_{\infty}^n}
\newcommand{\lif}{\ell_{\infty}}

\newcommand{\lrx}{\ell_{r,\infty}^{(x)}}
\newcommand{\lze}{\ell_{2,1}^{(a)}}

\newcommand{\Lrx}{{\cal L}_{r,\infty}^{(x)}}
\newcommand{\Lpa}{{\cal L}_{p,\infty}^{(a)}}
\newcommand{\cgi}{\gamma_{\infty}^{\it o}}
\newcommand{\Cgi}{\Gamma_{\infty}^{\it o}}

\newcommand{\tij}{T_{ij}}

\newcommand{\leb}{{\cal B}}
\newcommand{\B}{{\cal B}}
\newcommand{\com}{\B}

\renewcommand{\L}{{\cal L}}
\newcommand{\A}{{\cal A}}

\newcommand{\U}{{\cal U}}

\newcommand{\es}{{\cal S}}

\newcommand{\ef}{{\cal F}}

\renewcommand{\U}{{\cal U}}

\newcommand{\noo}{\left \|}
\newcommand{\rrm}{\right \|}

\newcommand{\bet}{\left |}
\newcommand{\rag}{\right |}

\newcommand{\intt}{\int\limits}
\newcommand{\summ}{\sum\limits}
\newcommand{\limm}{\lim\limits}

\newcommand{\pq}{\pi_{pq}}
\newcommand{\pr}{\pi_{r1}}
\newcommand{\opr}{\pi_{r1,cb}}
\newcommand{\opo}{\pi_{1,cb}}

\newcommand{\qed}{\hspace*{\fill}$\Box$\hz\pagebreak[1]}

\newcommand{\grc}{\gamma_{p,RC}}
\newcommand{\Grc}{\Gamma_{p,RC}}

\newcommand{\ioc}{\iota_{2,\infty}^k}

\documentstyle[11pt]{article}
\begin{document}

\begin{center}

{\bf \huge Orlicz property of operator spaces and eigenvalue
estimates}\hz

{\large Marius Junge}
\end{center}
\hz\hz
\begin{abstract}
\hspace{-1,9em} As is well known absolute convergence and
unconditional convergence for series are equivalent only in finite
dimensional
Banach spaces. Replacing the classical notion of absolutely summing
operators by the notion of 1 summing operators
\[ \summ_k \noo Tx_k \rrm \kl c \pl \noo \summ_k e_k \otimes x_k
\rrm_{\ell_1\otimes_{min}E}\]
in the category of operator spaces, it turns out that there are quite
different interesting
examples of 1 summing operator spaces. Moreover, the eigenvalues of a
composition
$TS$ decreases of order $n^{\frac{1}{q}}$ for all operators $S$
factorizing
completely through a commutative $C^*$-algebra if and only if the 1
summing norm
of the operator $T$ restricted to a $n$-dimensional subspace is not
larger than
$c n^{1-\frac{1}{q}}$, provided $q>2$. This notion of 1 summing
operators
is closely connected to the notion of minimal and maximal operator
spaces.
\end{abstract}

\section*{Introduction}
\newtheorem{heorem}{Theorem}
\newtheorem{propo}[heorem]{Proposition}
In Banach space theory the Orlicz property and its connection to
unconditional
and absolute convergence is well understood. For instance,
unconditional convergence and absolute convergence only coincide in
finite
dimensional spaces, but unconditional converging series are at least
square summable in the spaces $L_q$, $1\le q\le 2$. This was discovered
by Orlicz, hence the name Orlicz property. Furthermore,
this is best possible in arbitrary infinite dimensional Banach space by
Dvoretzky's theorem.
In the category of operator spaces there are several possibilities
to generalize the classical Orlicz property. We have choose a
definition
where only sequences are involved and which is motivated
by the theory of absolutely summing operators
introduced by Grothendieck. To be more precise, we recall
that an unconditional converging series $(x_k)_k \subset E$ in a Banach
space $E$ corresponds to an operator defined on $c_0$ with values in
$E$.
This is a consequence of the contraction principle.
\[ \noo \summ_k e_k \otimes x_k \rrm_{\ell_1 \otimes_{\eps} E}
\lel \sup_{|\alpha_k|\le 1} \noo \summ_k \alpha_k x_k \rrm \kl 4 \pl
 \sup_{\eps_k=\pm 1} \noo \summ_k \eps_k x_k \rrm \pl .\]
In order to involve the operator space structure of a an operator space
$E \subset B(H)$
we define an operator $T:E \nach F$ to be 1-summing if there exits a
constant
$c>0$ such that
\[ \summ_k \noo Tx_k\rrm \kl c \pl \noo \summ_k e_k \otimes x_k
\rrm_{\ell_1\otimes_{\min}E} \pl .\]
The best possible constant will be denoted by $\opo(T)$.
Here $\min$ denotes the minimal or spatial tensor product and $\ell_1$
is considered
as an operator space (for example by identification of the unit vectors
with the generators of a free group).
Anyhow, in the definition of absolutely summing operators
we simple replace the \underline{norm} by the \underline{cb-norm} of
the corresponding operator.
Obviously, in this definition only the Banach space structure of $F$ is
involved.
That's why this notion lives from the interplay
of operator space and Banach space theory.
For a more complete notion which is entirely
defined in the category of operator spaces and where matrices instead
of sequences
are considered we refer to the work of Pisier about completely
p-summing operators
and factorization problems. With this background 1-summing operators
turns to be the
weakest possible notion. The classical notion of absolutely
summing operators is included, by defining an operator spaces structure
of a Banach space via the embedding of $E$ in the commutative
$C^*$ algebra $C(B_{E^*})$. Following Paulsen we will denote this
operator
space by $\min(E)$.
In the first chapter we collect basic properties of 1 summing operators
and study the relation  between 1-summing operators and
$(1,C^*)$-summing operators defined early by Pisier for
$C^*$-algebras.
This is connected with Haagerup's characterization of injective von
Neumann algebras.

In this paper the framework of eigenvalue estimates
for operators factorizing completely through a commutative
$C^*$-algebra
is used to distinguish different operator spaces. The first part is
based on a
generalization of Maurey's inequality:

\begin{heorem} Let $2<q<\infty$. For an operator space $E$ a Banach
space $F$ and $T:E\nach F$ the following assertions are equivalent.
\begin{enumerate}
\item[i)] There exists a constant $c_1$ such that for all $n$
dimensional subspaces $G \subset E$ and $(x_k)_k \subset G$
\[ \summ_k\noo T x_k \rrm \kl c_1 \pl n^{1-\frac{1}{q}} \pl  \noo
\summ_k e_k \otimes x_k \rrm_{\ell_1 \otimes_{min} E} \pl. \]
\item[ii)] There exists a constant $c_2>0$ such that for all $\nen$ and
$x_1,..,x_n \in E$ one has
\[ \summ_1^n\noo T x_k \rrm \kl c_2\pl n^{1-\frac{1}{q}} \pl  \noo
\summ_1^n e_k \otimes x_k \rrm_{\ell_1^n \otimes_{min} E} \pl. \]
\item[iii)] There exists a constant $c_3$ such that for all operators
$S:F\nach E$ which factor completely
through a $C(K)$ space, i.e.  $S\p=\p PR$, $R:F \nach C(K)$ bounded and
$P: C(K)\nach E$
completely bounded one has
\[ \sup_{\nen} n^{\frac{1}{q}}|\lambda_n(ST)| \kl c_3\pl \noo P
\rrm_{cb} \pl \noo R: F\nach C(K)\rrm_{op} \pl.\]
\item[iv)]  There exists a constant $c_3$ such that for all operators
$S :F \nach E$ which factors completely through $B(H)$, i.e.
$S\p=\p PR$, $R:F \nach B(H)$ and $P: B(H)\nach E$  completely bounded
one has
\[ \sup_{\nen} n^{\frac{1}{q}}|\lambda_n(ST)| \kl c_4 \pl \noo P
\rrm_{cb} \pl \noo R: {\rm min}(F)\nach B(H)\rrm_{cb} \pl.\]
\end{enumerate}
where  $(\lambda_n(ST))_{\nen}$ denotes the sequence of eigenvalues in
non-increasing order according to their multiplicity.
\end{heorem}

As an application for identities we see that projection
constant of an $n$-dimensional subspace in a $(q,1)$-summing Banach
space
is at most $n^{\frac{1}{q}}$. This can already be deduced from
a corresponding theorem for identities on Banach spaces which was
proved in \cite{J}.
In the operator spaces setting we see that estimates for the growth
rate
of the 1-summing norm are useful to  measure
the 'distance' of an  operator space and its subspaces to
$\ell_{\infty}$ spaces.

Let us note that the theorem is not valid
for values $q<2$. For identities on Banach space this is not relevant,
all the properties are only satisfied for finite dimensional spaces.
In contrast to this operators spaces with 1-summing identity
are interesting spaces. For example the generators of the Cllifford
algebra spans such a space CL. This is probably not so surprising,
since this
example has proved to be relevant also for the closely connected notion
of
$(2,oh)$-summing operators, introduced by Pisier. Starting with CL we
construct
a scale of operator spaces with different growth rates for the
1-summing norm.
Further examples with small
1-summing norm were given
by randomly chosen $n$-dimensional
subspaces of the matrix algebra $M_N$, provided
$n\le N$. This was the starting point
to discover, independently of Paulsen and Pisier, the fact there are
only
few completely bounded operators between minimal and maximal operator
spaces.
Paulsen studied all possible operator space structures
on a given Banach space $E$ and realized that there is a minimal and
maximal one.
The minimal one is given by the commutative structure already defined
and the maximal by
the embedding $E \hookrightarrow  (\min(E^*))^*$, where $*$ denotes the
operators
space dual discovered  by Effros/Ruan and Blecher/Paulsen. This is
called the maximal
operator space $\max(E)$. Our approach is contained in the following
proposition which is a refinement of Paulsen/Pisier's result,
unfortunately
with a worse constant.

\begin{propo} Let $E$ be a maximal and $F$ be a minimal operator
space.
For an operator $T:F \nach E$ of rank at most $n$ one has
\[ \noo T \rrm_{cb} \kl \gamma_2^*(T) \kl 170 \noo T \otimes Id_{M_n}:
M_n(F) \nach M_n(E) \rrm \pl ,\]
where $\gamma_2^*$ is defined by trace duality with respect to Hilbert
space factorizing norm $\gamma_2$. In particular,
\[ \frac{\sqrt{n}}{170} \kl \noo Id:\min(E) \nach \max(E) \rrm_{cb}\pl
,\]
for all $n$ dimensional Banach space $E$.
\end{propo}

This is contained in the second part of this paper
where the study of 1-summing operators is continued.
This turns out to be quite fruitful in the context of dual operator
spaces.
For instance, maximal operator spaces are 1-summing if and only if they
are isomorphic to Hilbert spaces. Moreover, the 1-summing norm
of a $n$ dimensional subspace of $\max(\ell_r)$, $\max(\ell_{r'})$,
$\max(\es_r)$
$\max(\es_{r'})$ is less then $ 4 n^{\frac{1}{2}-\frac{1}{r}}$ for all
$2\le r\le \infty$.
Most of the techniques for maximal operator spaces carry over to duals
of exact operator
spaces by the key inequality of \cite{JP}. The lack of local
reflexivity
in operator spaces leads to the notion of exactness defined by Pisier
and motivated
by Kirchberg's work. In this context we will say that an operator space
is exact if all
its finite dimensional subspaces are uniformly cb-isomorphic to
subspaces of the matrix algebra's $M_N$. In the next theorem the
connection
between 1-summing operators and factorization properties is established
for duals of exact operator spaces.\hz

\begin{heorem} Let $1< p< 2$, $G$ an exact operator space, $E\subset
G^*$ and $F$ a minimal operator space.
For an operator $T:E \nach F$ the following are equivalent.
\begin{enumerate}
\item[i)] There exists a constant $c_1>0$ such that
\[ \summ_1^n \noo Tx_k  \rrm \kl c_1 \pl n^{1-\frac{1}{p}} \pl \noo
\summ_1^n e_i \otimes x_i \rrm_{\ell_1^n \otimes_{min} E} \pl. \]
\item[ii)] There is a constant $c_2$ such that for all completely
bounded operators $S:F \nach E$ one has
\[ \sup_k k^{\frac{1}{p}} \pl \bet \lambda_k(TS) \rag
\kl c_3 \pl \noo S:{\rm min}(F) \nach E \rrm_{cb} \pl .\]
\end{enumerate}
In the limit cases $p=1$ the eigenvalues are summable if and only if
the operator is 1-summing. In this case there exists a 1-summing
extension
$\hat{T}\p:\p G^* \nach \min(F^{**})$ which factors completely bounded
through $R\cap C$.
Furthermore, every completely bounded $S: \min{F} \nach E$ is
absolutely 2-summing
and hence the eigenvalues of a composition $TS$ are in $\ell_2$.
\end{heorem}

Let us note an application for an exact space $E\subset B(H)$ with
quotient map
$q:B(H)^*\nach E^*$. The Banach space $E^*$ is of cotype 2 and
$B(\ell_{\infty},E^*) \subset
CB(\ell_{\infty},E^*)$ if and only if there is a constant $c>0$ such
that for every sequence $(x_k)_1^n\subset E^*$ there is a sequence
$(\tilde{x}_k)_1^n \subset B(H)^*$
such that $q(\tilde{x}_k)=x_k$ and
\[ \ez \noo \summ_1^n \tilde{x}_k \eps_k \rrm_{B(H)^*} \kl c \pl  \ez
\noo \summ_1^n x_k \eps_k \rrm_{E^*} \pl .\]
In particular, a space $\max(X)$ is of operator cotype 2 if and only if
it satisfies
the conditions above if and only if it is a cotype 2 space satisfying
Grothendieck's theorem.
A non trivial example is the dual $A(D)^*$ of the disk algebra.

{\bf Acknowledgement:} I would like to thank Gilles Pisier for
stimulating
discussions, access to his preprints and the new concepts in operator
space
theory.

\setcounter{section}{0}
\section*{Preliminaries}
In what follows $c_0, c_1,$ .. always denote universal constants.
We use standard Banach space notation. In particular, the classical
spaces $\ell_q$ and $\ell_q^n$, $1\le q\le \infty$, $\nen$, are defined
in the usual way. We will also use the Lorentz spaces
$\ell_{p\infty}$.
This space consists of all sequences
$\si \in \lif$ such that
\[ \noo \si \rrm_{p\infty}  \pl := \pl  \sup_{\nen} \,
n^{\frac{1}{p}}\,\si_n^*
   \pl <\pl \infty . \]
Here $\si^*\,=\,(\si_n^*)_{\nen}$ denotes the
non-increasing rearrangement of $\si$.
The standard reference on operator ideals is the monograph of Pietsch
\cite{PIE}. The ideals of linear bounded operators, finite rank
operators,
integral operators are denoted by $\B$, $\ef$, ${\cal I}$. Given an
operator ideal
$(A,\alpha)$ the adjoint operator ideal $(A^*,\alpha^*)$ is defined by
the set of bounded
operators $T:Y\nach X$ such that
\[ \alpha^*(T) \pl:=\pl \sup\left \{ \bet tr(ST)\rag \mitt S\in
\ef(X,Y),\p \alpha(S) \kll 1\right \} \]
is finite. In particular, the ideal of integral operator is adjoint to
bounded operators with
\[ \iota_1(T) \pl := \pl  \noo \cdot \rrm^*(T) \pl .\]
We recall that an operator $T \in B(X,Y)$ factors through a Hilbert
space
($T\in \Gamma_2(X,Y)$)
if there are a Hilbert spaces $H$ and operators $S:X\nach H$, $R:H\nach
Y^{**}$
such that $\iota_{Y^*}T\lel RS$, where $\iota_{Y^*}:Y \nach Y^{**}$ is
the canonical
embedding of $Y$ into its bidual. The corresponding norm $\gamma_2(T)$
is defined as $\inf \{\noo S \rrm \noo R\rrm\}$, where the infimum is
taken over such
factorizations.
\hz

Let $1 \le q \le p \le \infty$ and $\nen$. For an operator $T \in
\leb(X,Y)$
the pq-summing norm of T with respect to $n$ vectors is defined by
\[ \pq^n(T) \pl := \pl
 \sup\left\{\p \kla \summ_1^n \noo Tx_k \rrm^p \mer^{1/p} \p \bet \pl
  \sup_{\noo x* \rrm_{X^*}\le1} \kla \summ_1^n \bet\langle
  x_k,x^*\rangle \rag^q \mer^{1/q}
 \pl \le \pl 1\right.\p \right\} \pl .\]
An operator is said to be absolutely pq-summing
$(T \in \Pi_{pq}(X,Y))$ if
\[ \pq(T) \pl := \pl \sup_n \pq^n(T) \pl < \pl \infty \pl . \]
Then $(\Pi_{pq},\pq)$ is a maximal and injective Banach ideal (in the
sense
of Pietsch). As usual we abbreviate  $(\Pi_q,\pi_q) :=
(\Pi_{qq},\pi_{qq})$.
For further information about absolutely pq-summing operators we refer
to
the monograph of Tomczak-Jaegermann \cite{TOJ}.

The definition of some s-numbers of an  operator $T\in\com (E,F)$ is
needed.
The $n$-th $approximation$ $number$ is defined by

\[ a_n(T) \pl :=\pl \inf\{\, \noo T-S \rrm \, | \,rank(S)\,< \, n \,\}
\pla ,\]

whereas the $n$-th $Weyl\,number$ is given by

\[ x_n(T) \pl :=\pl \sup\{\, a_n(Tu)\, |\, u \in \B(\ell_2,E) \,
\mbox{with}
 \, \noo u \rrm \,\le \, 1\,\} \pla .\]
Let $s \in \{a,x\}$. By ${\cal L}_{pq}^{(s)}$ we denote the ideal of
operators $T$
such that $(s_n(T))_{\nen}\in \ell_{pq}$ with the associated quasi-norm
$\ell_{pq}^{(s)}(T) \p:=\p \noo (s_n(T))_{\nen} \rrm_{\ell_{pq}}$.
If $H$ is a Hilbert space the spaces ${\cal S}_{pq}(H) \p=\p {\cal
L}_{pq}^{(a)}$
are normable. Indeed all $s$-numbers coincide for operators
on Hilbert spaces. If $p\lell q$ we will briefly write ${\cal
S}_p(H)$.
This includes ${\cal S}_2(H)$ the set of Hilbert-Schmidt operators.
\hz

By Ruan characterization theorem there are two possibilities to
introduce operator spaces. Either as subspaces of $\B(H)$, where $H$ is
a
Hilbert space or as a Banach space $E$ together with a sequence of
norms
on the spaces of $n\times n$ matrices $M_n(E)$ with values in $E$.
To guarantee that such a sequence of norms is induced by
an embedding into some $B(H)$ the following axioms are required.
\begin{enumerate}
\item[i)] If $O=(O_{ij})$, $P=(P_{ij})$ are scalar $n\times n$ matrices
and $x =(x_{ij})$ in $M_n(E)$
one has
\[ \noo (\summ_{kl} O_{ik}x_{kl}P_{lj})_{ij} \rrm_{M_n(E)} \kl \noo
O\rrm \pl \noo x \rrm_{M_n(E)}\pl \noo P\rrm\pl .\]
\item[ii)] If a matrix $B \lell \kla {x \atop 0}{0\atop y}\mer$
consists of two
disjoint blocs one has
\[ \noo B \rrm \lel \max\{ \noo x \rrm, \p \noo y \rrm \} \pl .\]
\end{enumerate}
A major step for the development of operator space theory
is the right definition of an operator space dual. Indeed, the norm
of a matrix $(x^*_{ij})\subset E^*$ is given by
\[ \noo x^*_{ij} \rrm_{M_n(E^*)} \lel \noo (x^*_{ij}) :E\nach M_n
\rrm_{cb}
\lel \sup \left\{ \noo \langle x^*_{ij}, x_{kl}\rangle \rrm_{M_{n^2}}
\mitt \noo x_{ij} \rrm_{M_n(E)} \kll 1 \right\} \pl .\]
For further information on this and operator space theory we
refer to the paper of Blecher and Paulsen, \cite{BPT}.

\section{The notion of 1-summing operators on operator spaces}

Given two Banach spaces $X$ and $Y$ a matrix structure corresponding
to operator spaces is defined on $\B(X,Y)$ in the following way.
The norm of a matrix $(T_{ij})\subset \B(X,Y)$ is induced
by considering this matrix as element in $\B(\lzn(X),\lzn(Y))$
\[ \noo T_{ij} \rrm_n \pl := \pl \sup\left\{ \kla \summ_{i=1}^n
\noo \summ_{j=1}^n \tij (x_j)\rrm^2 \mer^{\frac{1}{2}} \pl\left |\pl
\summ_1^n \noo x_j \rrm^2 \pl
 \le \pl 1 \pl \right. \right\}\pl . \]
Following \cite{PCB}
an operator $u \in \B(E,F)$, where $E \subset \B(X_1,Y_1)$ and
$F \subset \B(X_2,Y_2)$ is said to be completely bounded if there is a
constant
$c>0$ such that for $(\tij) \subset E$
\[ \noo u(\tij) \rrm_n \pl \le \pl c \pl \noo \tij \rrm_n\pl . \]
The infimum over all such constants is denoted by $\noo u \rrm_{cb}$.
As usual $\lin$ will be considered as a subspace
of $\B(\lzn)$. The matrix norm induced by this embedding corresponds to
the $\eps$ tensor product. In analogy to the classical theory of
absolutely r1-summing operators we
define the r1-summing norm (with $n$ vectors) for an
operator $T \in B(E,F)$, where $F$ is a Banach space and
$E \subset B(X,Y)$ as follows
\for
\opr^n(T) &:=& \sup\left \{ \kla \summ_1^n \noo Tu(e_k) \rrm^r
\mer^{\frac{1}{r}}
\pl \left | { \atop } \noo u\p:\p\lin \nach E \rrm_{cb} \pl \le 1
\pl\right.\right\} \pll \mbox{and}\\
\opr(T) &:=& \sup_{\nen} \opr^n(T) \pl .
\mel
An operator $T$ is said to be $r1-summing$ if $\opr(T)$ is finite.
The notion of absolutely r1-summing operators is included in this
definition if we consider $E$ to be embedded into $C(B_{E^*})\subset
\B(\ell_2(B_{E^*}),\ell_2(B_{E^*}))$.
A basic tool for the notion of r1-summing operators is a description
of the cb norm for operators acting on $\lif$. This is well-known but
since
it is crucial for the following we give a proof.

\begin{lemma}\label{cund} Let $E \subset \B(X,Y)$ and $u \in
\B(\lin,E)$ with $x_k=u(e_k)$.
Then we have
\[ \noo u \rrm_{cb} \pl = \pl \sup\left\{ \summ_1^n
\sigma_1(vx_kw)\pl\mitt\pl v\in \B(Y,\ell_2),\pl w\in \B(\ell_2,X) \pl
\mbox{and} \pl\pi_2(v),\pi_2(w^*)\pl\le\pl 1 \pl\right\}\pl , \]
where $\sigma_1$ denotes the trace class norm.
\end{lemma}

{\bf Proof:} Clearly, the supremum on the right hand remains unchanged
if
we replace all operators $v \in \B(Y,\ell_2)$, $w\in \B(\ell_2,X)$
by the supremum over $m\in \nz$ and $v\in \B(Y,\ell_2^m)$,
$u\in \B(\ell_2^m,X)$. By a well known characterization of 2-summing
operators, see \cite{PIL}, every operator $v\in \B(Y,\ell_2^m)$ can
be written in the form $v=Oz$ with
\[ \kla \summ_1^N \noo z^*(e_k) \rrm^2 \mer^{\frac{1}{2}}\noo O\p: \p
\ell_2^N \nach \ell_2^m
\rrm \pl \le \pl (1+\eps) \pl \pi_2(v)\pl , \]
for $\eps>0$ arbitrary. Hence we get
\for\lefteqn{
\sup\left\{ \summ_1^n \sigma_1(vx_kw)\pl\mitt\pl
v\in \B(Y,\ell_2),\pl w\in \B(\ell_2,X) \pl \mbox{and}
\pl\pi_2(v),\pi_2(w^*)
\pl\le\pl 1 \pl \right\}\pl = }\\
&=& \sup_{N\in \nz} \sup\left\{ \summ_1^n tr(A^k vx_kw) \mitt \pl
\noo A^k \p :\p\ell_2^N \nach \ell_2^N \rrm \pl \le \pl 1,
\pl \summ_1^N \noo w(e_i)
 \rrm^2,\pl \summ_1^N \noo v^*(e_j) \rrm^2 \pl \le \pl 1  \right\}\\
&=& \sup_{N \in \nz} \sup \left \{ \summ_{k=1}^n \summ_{i=1}^N
<v^*(e_i), \summ_{j=1}^N A^k_{ji}x_k(w(e_j))>  \mitt \pl
\noo A^k\p:\p \ell_2^N \nach \ell_2^N \rrm \pl \le \pl 1,\pl
\summ_1^N \noo w(e_i) \rrm^2, \right.\\
& &\pla \pla \pla \pla\pla \pla\pla\pll\pll \left.\summ_1^N \noo
v^*(e_j) \rrm^2 \pl \le \pl 1  \right\}\\
&=&\sup_{N \in \nz} \sup\left\{ \noo \kla u\kla \summ_1^n e_k \otimes
A^k_{ji}\mer \mer_{ij} \rrm_N \mitt\pll  \sup_k \noo (A^k)^t \rrm \pl
\le \pl 1 \right\}\\
&=& \noo u \rrm_{cb}\pl .\\[-1.3cm]
\mel\hfill $\Box$\hz
\begin{rem} \label{rcoh}
{\rm If $E \subset X^* \cong B(X,\cz)$ or $E \subset Y\cong B(\cz,Y)$
the formula above
reduces to
\[ \noo u :\lin \nach E \rrm_{cb} \pl =\pl \pi_2(u) \pl .\]
Therefore the 1-summing norm of an operator $T \in B(E,F)$
coincides with the absolutely 2-summing norm
\[ \opo(T) \pl =\pl \pi_2(T) \pl .\]
If the space $E$ has Cotype 2 (or is $(2,1)$ mixing, see \cite{PIE})
every
absolutely 2-summing operator is absolutely 1-summing and therefore all
these notions coincide. The most canonical examples are given by the
row space
$R = \B(\cz,\ell_2)$ and  the column space $C = \B(\ell_2,\cz)$.
In this cases it is a consequence of the "little Grothendieck
inequality", see \cite{TOJ},
\[ \pi_1(T) \kl \pl \frac{2}{\sqrt{\pi}}\pll \pi_2(T) \pl =\pl
\frac{2}{\sqrt{\pi}}\pll \opo(T)\pl .  \]
By interpolation the same remains true for the operator Hilbert space
$OH$.}
\end{rem}

{\bf Proof:} Let $E \subset Y\cong \B(\cz,Y)$ and $u \in \B(\lin,E)$.
Trace duality for the absolutely 2 summing operators implies
\for
\noo u \rrm_{cb} &=& \sup_{\pi_2(v),\pi_2(w^*)\le 1} \summ_1^n
\si_1(v(e_1 \otimes y_i)w)\kl
\sup_{\pi_2(v),\pi_2(w^*)\le 1} \summ_1^n \noo v(y_i) \rrm  \pl \noo w
\rrm\\
&=&\sup_{\pi_2(v)\le 1} \iota_1(vu)\kl
\sup_{\pi_2(v), \noo w \rrm \le 1} \bet tr(vuw)\rag\kl
\sup_{\noo w \rrm  \le 1} \pi_2(uw) \lel \pi_2(u) \pl .
\mel
The argument for $E \subset \B(X,\cz)$ is similar. For $T \in B(E,F)$
we use Pietsch' factorization theorem, again trace duality
and the fact that absolutely 1-summing operators on $\ell_{\infty}$
are integral
\for
\opo(T) &=& \sup \left \{ \pi_1(Tu) \mitt \pi_2(u: \ell_{\infty}^n
\nach E) \le 1\right\}\\
&=& \sup\left \{ \iota_1(Tu) \mitt \pi_2(u: F \nach E) \le 1\right \}
\lel \pi_2(T)\pl . \\[-1.5cm]
\mel \qed \hz

Nowadays it can be considered as a standard application of the
Hahn-Banach
separation theorem to deduce a factorization theorem for
1-summing operators. We refer to \cite{PSP} for the required
modification
in the infinite dimensional case.
\begin{prop}\label{fac} Let $X$, $Y$, $F$ be Banach spaces,
$E \subset \B(X,Y)$ and  $u \in \B(E,F)$.
\begin{enumerate}
\item Let us assume that $X$ and $Y$ finite dimensional, of dimension
$n$ and $m \in \nz$, say. The operator $T$ is 1-summing if and only if
there exists
a constant $C>0$ and a probability measure $\mu$ on the compact space
$K\p:=\p {\rm B}_{\Pi_2^d(\ell_2^n,X)} \times {\rm
B}_{\Pi_2(Y,\ell_2^m)}$
such that
\[ \noo Tx \rrm \pl \le \pl C \intt_K \sigma_1(vxu) \pl d\mu(u,v) \pl .
\]
\item $T$ is 1-summing  if and only if there exists a constant $C>0$
and
an ultrafilter $\U$ over an index set $\A$ together with finite
sequences
$(\lambda^{\alpha}_i)_{i \in I^{\alpha}}$,
$(u_i^{\alpha}, v_i^{\alpha})_{i \in I^{\alpha}} \subset
\B_{\Pi_2^d(\ell_2,X)}\times B_{\Pi_2(Y,\ell_2)}$
such that
\[ \noo Tx \rrm \pl \le \pl C \limm_{\alpha \in \U} \summ_{i \in
I^{\alpha}}  \sigma_1(v_i^{\alpha}xu_i^{\alpha}) \pl . \]
\end{enumerate}
In both cases  C can be chosen to be $\opo(T)$. In particular if
$E\subset B(H)$
is an operator space and $F$ carries its minimal (commmutative)
operator spaces
structure then every 1-summing operator is completely 1-summing in the
sense of Pisier, \cite{PSP}.
\end{prop}



In the next proposition we list the relations between the notion of
r1-summing operators and $(r1,C^*)$-summing operators defined on
$C^*$-algebra's by Pisier. More generally, let us recall that
an element $z\in \B(X,\ov{X^*})$, $\ov{X^*}$ the anti dual, is said to
be
positive if $\langle z(x),x \rangle\gll 0$ for all $x \in X$. An
operator
$u:\lin \nach \B(X,\ov{X^*})$ is positive, if $u$ maps positive
sequences
into positive elements.

\begin{prop} Let $X$ be a Banach space.
\begin{enumerate}
\item An operator $u:\lin \nach \B(X,\ov{X^*})$ is completely bounded
if and only if $u$ is decomposable into positive operators and
\[ \noo u \rrm_{cb} \kl \inf\left\{ \summ_j \bet \lambda_j\rag \noo u_j
\rrm\mitt u \p=\p \summ_j \lambda_j u_j\pl ,\pl u_j \p positive\right\}
\kl 4 \pl \noo u \rrm_{cb} \pl .\]
Therefore an operator $T :\B(X,\ov{X^*}) \nach F$ is r1-summing if and
only if
\[\kla \summ_1^n \noo Tz_k \rrm^r \mer^{ \frac{1}{r} }  \kl C \pl \noo
\summ_1^n z_k \rrm\]
for all finite sequences of positive elements $(z_k)_1^n \subset
\B(X,\ov{X}^*)$.
The corresponding constants are equivalent up to a factor 4. Given an
operator
$v:X^* \nach G$ then  operator $T:=v \otimes \bar{v}:\B(X ,\ov{X^*})
\nach G\otimes_{\eps}\ov{G}$
is 1-summing if and only if $v$ is absolutely 2-summing.
\item If $E$ is a subspace of a $C^*$-algebra and $T \in \B(E,F)$
is a r1-summing operator then it is $(r1,C^*)$-summing, i. e.
for all $(x_k)_k \subset C^*$
\[ \kla \summ^n_1  \p \noo u(x_k) \rrm^r\mer^{\frac{1}{r}} \pl \le \pl
4 \pl \opr(T) \pl \noo
\summ_1^n \kla \frac{x^*_k x^{ }_k + x^{ }_k x_k^*}{2}
\mer^{\frac{1}{2}} \rrm_{C^*} \pl. \]
Conversely, if $E$ is a von Neumann algebra, $E$ is injective if and
only if
every $(1,C^*)$-summing operator is 1-summing and satisfies
\[ \opo(T) \kl c \pl \pi_{1,C^*}(T) \pl ,\]
where $c$ is a constant depending on $E$ {\rm (}$\pi_{1,C^*}$ denotes
the best constant in the inequality above for $r=1${\rm )}.
In this case also $\opr(T) \kll c \p \pi_{r1,C^*}(T)$ for all $1\le
r<\infty$.

\item If $T :E\nach F$ is a 1-summing operator defined on an operator
space
$E \subset \B(H)$ it is $(2,oh)$, $(2,R)$ and $(2,C)$-summing. This
means
\for
 \kla \summ^n_1  \p \noo T(x_k) \rrm^2  \mer^{\frac{1}{2}}&\le& \pl
 \opo(T) \pl
\noo  \summ_1^n  x_k \otimes \overline{x_k} \rrm_{E \otimes_{min}
\overline{E}}^{\frac{1}{2}}\pl ,\pla\\
\kla \summ^n_1  \p \noo T(x_k) \rrm^2  \mer^{\frac{1}{2}}&\le& \pl
\opo(T) \pl
\noo  \kla \summ_1^n  x^{ }_k x^*_k \mer^{\frac{1}{2}}\rrm_{\com(H)}\pl
, \pla\\
\mbox{and}\pla
\kla \summ^n_1  \p \noo T(x_k) \rrm^2  \mer^{\frac{1}{2}}&\le& \pl
\opo(T) \pl
\noo  \kla \summ_1^n  x^*_k x^{ }_k
\mer^{\frac{1}{2}}\rrm_{\com(H)}\pl. \pla \\
\mel
\end{enumerate}
\end{prop}\hz

{\bf Proof:} For the following let us denote by $\pi_{r1}^+(T)$ the
best constant $C$ satisfying
\[ \kla \summ_1^n \noo T(z_k)  \rrm_{F}^r \mer^{\frac{1}{r}} \pl \le
\pl C  \pl  \noo \summ_1^n z_k \rrm_{\B(X,\ov{X^*})}  \]
for all positive elements $(z_k)_1^n$. Then we have trivially
\[ \kla \summ_1^n \noo T(u(e_k))  \rrm_{F}^r \mer^{\frac{1}{r}} \pl \le
\pl \pi_{r1}^+(T) \pl  \noo u :\lin \nach \B(X,\ov{X^*}) \rrm_{dec} \]
where
\[ \noo u \rrm_{dec} \pl :=\pl \inf\left\{ \summ_j \bet \lambda_j\rag
\noo u_j \rrm_{op}\mitt u \p=\p \summ_j \lambda_j u_j\pl ,\pl u_j\p
positive\right\} \pl .\]
We will first show that for a positive operator $u$
\[ \noo u \rrm_{cb} \lel \noo \summ_1^n u(e_k) \rrm \lel \noo u
\rrm_{op} \pl .\]
For this we can assume that $z_k \lell u(e_k)$
are positive elements in $\B(X,\ov{X^*})$. Let us note that positive
elements
are automatically  $\Gamma_2$ operators. On the tensor product
$\ell_2\otimes X$
we use the norm induced by the absolutely $2$ summing norm of the
corresponding
operator from $X^*$ with values in $\ell_2$. With this norm each
element
$x_k$ defines a positive, possibly degenerated, scalar product
\[ \phi_k \p:\kla\ell_2\otimes X\mer \times \kla \ell_2\otimes X\mer
\nach \cz\pl \quad\mbox{with}\quad
 \phi_k(v,w) \p:=\p tr(\ov{v^*}z_kw) \pl.\]
From Lemma \ref{cund}, H\"older's and the Cauchy-Schwartz inequality we
deduce
\for
\noo u \rrm_{cb} &=& \sup\left \{ \summ_1^n tr(\ov{A^k}\ov{v^*}z_kw)
\pl\mitt \pi_2(v^*), \pi_2(w^*), \noo A^k \rrm\pl \le \pl 1 \right\}\\
 &=&  \sup\left \{ \summ_1^n \phi_k(vA^k,w) \pl\mitt
 \pi_2(v^*), \pi_2(w^*), \noo A^k \rrm\pl \le \pl 1\right\}\\
 &\le&  \sup\left \{ \kla \summ_1^n \phi_k(vA^k,vA^k)
 \mer^{\frac{1}{2}}\p
\kla  \summ_1^n \phi_k(w,w) \mer^{\frac{1}{2}}
 \pl\mitt  \pi_2(v^*), \pi_2(w^*), \noo A^k \rrm\pl \le \pl 1\right\}\\
 &\le& \sup\left\{  \summ_1^n \si_1(\ov{v^*}z_kv) \pl \mitt
 \pi_2(v^*)\pl \le \pl 1 \right\}\lel
 \sup\left \{ \summ_1^n tr(\ov{v^*}z_kv) \pl \mitt \pi_2(v^*)\pl \le
 \pl 1 \right\}\\
 &\le& \gamma_2(\summ_1^n z_k) \kl \noo \summ_1^n z_k \rrm  \lel \noo
 u(1,..,1) \rrm \kl \noo u \rrm_{op} \pl .
\mel
Where we used that for a positive element $z_k$ the composition
$\ov{v^*}z_kv$ actually
defines a positive operator on $\ell_2$ and that for the positive
element $\summ z_k$
the $\gamma_2$-norm and the operator norm coincide.
(If $X \p=\p H$ the whole statement can be deduced from
\cite[theorem 2.4., proposition 3.5.]{PAU}.) In particular we obtain
\[ \noo u \rrm_{cb} \kl \noo u \rrm_{dec} \quad \mbox{and} \quad
\pi_{r1}^+(T) \kl \opr(T) \pl .\]
{\bf\it 1:} Let $u: \lin \nach \B(X,\ov{X^*})$ be a completely bounded
operator.
By Pisier's version \cite{PCB}
of the Haagerup/Wittstock factorization theorem, there exists
a $*$-representation $\pi: \B(\ell_2^n) \nach \B(H)$ and operators
$V,W: H\nach X$ such that
\[ u(\alpha) \lel \ov{W^*}\p \pi(D_{\alpha})\p V \quad \mbox{and} \quad
\noo V\rrm \lel \noo W \rrm \kl \sqrt{\noo u \rrm_{cb}} \pl ,\]
where $D_{\alpha}$ denotes the diagonal operator with entries $\alpha$.
It is standard
to see that the operators
\[ u^k(\alpha) \pl :=\pl \frac{1}{4}\ov{(V+i^k W)}\p  \pi(D_{\alpha})
\p (V+i^kW) \quad k=0,..,3 \]
are positive and of norm less than $\noo u \rrm_{cb}$. But $u \lel
u^0-u^2 + i(u^1-u^3)$
implies $\noo u \rrm_{dec} \kl 4 \noo u \rrm_{cb}$.
The second statement about operators $T$ of the form $v\otimes \ov{v}$
is a simple consequence of the observation that elementary tensors
$z_i \lell x^*_i \otimes \ov{x^*_i}$ are clearly positive. For the
reverse implication one simply uses Pietsch factorization theorem
for absolutely 2 summing operators.

{\bf\it 2:} Clearly we have $\pi_{r1}^+(T) \kl \pi_{r1,C^*}(T)$. For
the
converse we only have to note that every element $x$ in a $C^*$ algebra
admits a decomposition $x\lel x^1-x^2+i(x^3-x^4)$ in positive elements
such that
\[ x^k \kl  \kla \frac{x^* x + xx^*}{2} \mer^{\frac{1}{2}} \pl .\]
Hence we get $\pi_{r1,C^*}(T) \kl 4 \pi_{r1}^+(T)$. If $E$ is a von
Neumann
algebra we see that the existence of a constant $c_1>0$
\[ \opo(T) \kl c_1 \pl \pi_{1,C^*}(T) \pl ,\]
for all operators $T:E \nach \lin$ is equivalent with the existence of
a constant $c_2$
\[ \iota^o(T) \lel  \opo(T) \kl c_2\pl \pi_1^+(T) \pl ,\]
where $\iota^o$ is the operator integral norm.
Hence
trace duality implies that the condition above is equivalent to
\[ \noo u \rrm_{dec} \kl c_2 \noo u \rrm_{cb} \pl \]
for all $u:\lin \nach E$. By Haagerup's theorem, see \cite{HA}, this
holds if and only if $E$ is
injective. Together with the proof of
1.we see that for an injective von Neumann algebra the notion of
r1-summing
and $(r1,C^*)$-summing coincide.

{\bf\it 3:} This is an easy variant of Kwapien's argument. By the
remark \ref{rcoh} we deduce that
for all diagonal operator $\ds \p: \lin \nach \lzn$ and $G_n \in
\{R_n,C_n,OH_n\}$
\[ \noo \ds \p :\p\lin \nach G_n \rrm_{cb} \pl =\pl \pi_2(\ds) \pl =\pl
\noo \si \rrm_2 \pl .\]
Let us denote by $(e_k)_1^n$ the sequence of unit vectors of $G_n$.
Then we get for
all $w \in \com(G_n,E)$
\for
\kla \summ_1^n \noo Tw(e_k) \rrm^2 \mer^{\frac{1}{2}} &=&
\sup_{\noo \si \rrm_2 \p\le\p1} \summ_1^n \noo Tw\ds(e_k) \rrm \\
 &\le& \opo(T) \pl \sup_{\noo \si \rrm_2 \p\le\p1} \noo w\ds \rrm_{cb}
 \kl
 \opo(T) \pl \noo w \rrm_{cb} \pl.
\mel
The assertion is proved by identifying the complete bounded norm of w
with
the corresponding expressions on the right hand side in 2.. For
$G_n \p=\p OH_n$ this was done in \cite{PLT}. For the two other cases
we refer to
\cite{BPT}.\hfill $\Box$


\begin{rem}\label{ocot} {\rm For an operator space $E\subset \com(H)$
which is of operator
cotype 2 the a priori different notions of summability coincide.
Indeed, using the same
arguments as in the commutative theory, see {\rm \cite{PIL}}, one can
deduce that every operator
$S \in \com(\lin,E)$ factors through $OH_n$ with $\gamma_{oh}(S)\p\le
c(E)\p\noo S\rrm$. For notation and information see {\rm \cite{PLT}}.
A use of "little Grothendieck" inequality implies
\[ \pi_1(T) \pl \le \pl c_0\pl c(E)\pl \pi_{2,oh}(T)\pl . \]
For all (2,oh)-summing operator $T \in \com(E,\ell_2)$. Finally the
factorization properties of (2,oh)-summing operators imply for all
operators $T\in \com(E,F)$
\[ \frac{1}{c_0 c(E)}\pl \pi_1(T) \pl \le \pl \pi_{2,oh}(T)
\pl \le \pl \opo(T) \pl \le \pl \pi_1(T)\pl . \]
}
\end{rem}

The proof of the first theorem in the introduction
is based on a similar statement for the absolutely-summing norm
of operators defined on $C(K)$ spaces.

\begin{prop} \label{connect} Let $2 < r < \infty$, $K$ a compact
Haussdorf space, $F$ a Banach space
and $T:C(K) \nach F$. If there exists a constant $C>0$ such that
\[ \summ_1^n \noo Tx_k \rrm  \kl C \pl n^{1-\frac{1}{r}} \pl \sup_{t\in
K} \summ_1^n \bet x_k(t) \rag \]
for all elements $(x_k)_1^n \subset C(K)$, then we have
\[ \lrx(T) \pl \le \pl c_0 \p \kla \frac{1}{2}-\frac{1}{r} \mer^{-1}
\pl C \pl , \]
where $c_0$ is an absolute constant. If $F$ and $C(K)$ are complex
Banach spaces one has for  every $S:F \nach C(K)$
\[ \sup_{\ken} k^{1/r} \la_k(TS) \pl \le \pl c^2_0 \p \kla
\frac{1}{2}-\frac{1}{r} \mer^{-1} \pl C \pl \noo S \rrm \pl . \]
\end{prop}

{\bf Proof:} First we show
\[ \noo \pl (\noo Tx_k \rrm_F) \pl \rrm_{r,\infty} \pl \le \pl C
\pl \sup_{t \in K} \summ_k \bet x_k(t)\rag \pl  \]
for all $(x_k)_1^n \subset C(K)$. Indeed we can assume $\noo Tx_j\rrm$
non increasing. For fixed $1\kll k \kll n$ we get
\for
k \p \noo Tx_k \rrm &\le&
\summ_1^k \noo Tx_l \rrm \kl C \pl k^{1-\frac{1}{r}}
\pl
\sup_{t \in K} \summ_1^k \bet x_j(t)\rag
\mel
Dividing by $k^{1-\frac{1}{r}}$ and taking the supremum over all $1
\kll k \kll n$ yields the
estimate. Now we choose $2<q<r$ with
$\frac{1}{2}+\frac{1}{r}=\frac{2}{q}$.
For $(x_k)_1^n \subset C(K)$ we obtain
\for
\kla \summ_1^n \noo Tx_k \rrm^q \mer^{1/q} &\le&
\kla \summ_1^n k^{-q/r} \mer^{1/q} \pl (\noo \pl\noo Tx_k\rrm )\pl
\rrm_{r,\infty}\\
&\le& \kla \frac{1}{q}-\frac{1}{r} \mer^{-1/q} \pl n^{1/q-1/r} \pl
c_2 \pl \sup_{t \in K} \summ_1^n \bet x_k(t)\rag \pl .
\mel
Therefore we have
\[ \pi_{q1}^n(T) \pl \le \pl C
\pl \kla \frac{1}{q}-\frac{1}{r} \mer^{-1/q} \pl n^{1/q-1/r} \pl . \]
Using Maurey's theorem, see \cite[theorem 21.7]{TOJ}\p, this implies
with our
choice of q
\for
\pi_{q2}^n(T) &\le& C \pl c_0
\pl  \kla \frac{1}{2}-\frac{1}{q} \mer^{1/q-1} \pl
\pl \kla \frac{1}{q}-\frac{1}{r} \mer^{-1/q} \pl n^{1/q-1/r} \\
&\le& C \pl 2\p c_0 \pl \kla \frac{1}{2}-\frac{1}{r} \mer^{-1} \pl
n^{1/q-1/r} \pl .
\mel
Now let $u \in \leb(\ell_2,C(K))$. By a Lemma, probably due to Lewis,
see
\cite[Lemma 2.7.1]{PIE}, one can find for all $\nen$ an orthonormal
family $(o_k)_1^n$ in $\ell_2$ with
\[ a_k(Tu) \pl \le \pl 2\pl \noo Tu(o_k) \rrm \quad \mbox{for all}
\quad k=1,..,n\pl .\]
Hence we deduce
\for
n^{1/q} a_n(Tu) &\le& 2 \pl \kla \summ_1^n \noo Tu(o_k) \rrm^q
\mer^{1/q}\kl
2\pl \pi_{q2}^n(T) \pl \sup_{t \in K} \kla \summ_1^n \bet u(o_k)(t)
\rag^2 \mer^{1/2}\\
&\le& 4\pl C \pl c_0 \pl \kla \frac{1}{2}-\frac{1}{r} \mer^{-1} \pl
n^{1/q-1/r} \pl
\sup_{\noo \al \rrm_2 \le 1} \noo u\kla \summ_1^n a_k\p o_k\mer
\rrm_{C(K)} \\
&\le& 4\pl C \pl c_0 \pl \kla \frac{1}{2}-\frac{1}{r} \mer^{-1}\pl
n^{1/q-1/r}\pl \noo u \rrm\pl .
\mel
Dividing by the factor $n^{1/q-1/r}$ and taking the supremum over
$\nen$
yields
\[ \sup_{\nen} n^{1/r} \p a_n(Tu) \pl \le \pl 4\p c_0 \pl \kla
\frac{1}{2}-\frac{1}{r} \mer^{-1} \pl C \noo u \rrm \pl . \]
Now taking the supremum over all u with norm less than 1 the desired
estimate for
the Weyl numbers is proved. For the estimates of the eigenvalues
we use the fact that the ideal  $\Lrx$ is of eigenvalue type
$\ell_{r,\infty}$, \cite[3.6.5]{PIE}.
\qed

\begin{rem} {\rm In fact all these conditions are equivalent as far as
$2<r<\infty$.
If $1<r<2$ let us consider the embedding $I :\ell_1\nach  C[0,2\pi]$
given by the
Rademacher functions $r_j(t)= sign \sin(2^jt)$ and the corresponding
projection
$P: C[0,2\pi] \nach \ell_2$. By Kintchine's inequality $P$ is
r1-summing
for all $r>1$. On the other hand if we compose with a continuous
diagonal operator $D_{\tau}: \ell_2 \nach \ell_1$
we see that the best possible eigenvalue behaviour for r1-summing
operators
is actually $(\lambda_k(PD_{\tau}))_{\ken} \in \ell_2$. For $r=2$ a
more complicated example was
constructed by \cite{KOE}. This shows that the assumption $r>2$ is
really
necessary. }
\end{rem}

\begin{rem} {\rm Since for an operator $A \in \B(\lin,\ell_{\infty}^m)$
the
operator norm coincides with completely bounded norm we have for
$1 \le r \le \infty$
\[ \opr^n(u) \pl = \pl \sup\left\{ \pr^n(uw)\pl \mitt \pl \noo w\p:\p
\ell_{\infty}^m \nach E\rrm_{cb} \pl \le \pl 1 \right\} \pl . \]
Therefore the results of {\rm \cite{DJ}} can be applied to deduce for
each
operator u of rank at most $n$
\[\opr(u) \pl \le \pl c_0^{\frac{r'}{r}} \left \{\begin{array}{l
@{\quad} l}
 \kla \frac{1}{r}-\frac{1}{2} \mer^{-\frac{r'}{2r}} \pl
 \opr^{[n^{r'/2}]}(u) & for \pl 1 < r <2\\[+0.2cm]
 \pi_{21,op}^{[n(1+\ln n)]} (u) & for \pl r=2\\[+0.2cm]
 \kla \frac{1}{2}-\frac{1}{r} \mer^{\frac{1}{r}} \pl \opr^n(u)
	    & for \pl 2<r<\infty \pl ,
\end{array}\right.\]
where $r'$ is the conjugate index to $r$.}
\end{rem}

An operator $u \in \B(F,E)$, $E \subset \B(X,Y)$ is said to be
$completely$
$\infty-factorable$ ($u \in \Cgi(F,E)$) if there is a factorization
$u\p=\p
SR$, where $R\in CB(F,B(H))$, $S \in CB(B(H),E)$, H a Hilbert space.
The $\cgi$-norm of u is defined
as $\inf\{\noo S \rrm_{cb}\p\noo R\rrm_{cb}\}$ where the infimum is
taken over
all such factorizations. As in the commutative case this turns out to
be a
norm. Now we can prove the first theorem of the
introduction.

\begin{theorem} \label{opeigen} Let $2<r<\infty$, $X$, $Y$, $F$ Banach
spaces
and $E \subset B(X,Y)$. For an operator $T:E \nach F$
the following assertions are equivalent.
\begin{enumerate}
\item[i)] There is a constant $c_1$ such that for all $\nen$
\[ \opo^n(T) \pl \le \pl c_1 \pl n^{1-\frac{1}{r}} \pl . \]
\item[ii)] There is a constant $c_2$ such that for all operators $R\in
\B(F,C(K))$,
$S\in CB(C(K),E)$, $K$ a compact Haussdorf space
\[ \sup_{\ken} k^{1/r} \bet \la_k(TSR) \rag \pl \le \pl c_2 \pl \noo R
\rrm\pl
\noo S\rrm_{cb} \pl . \]
\item[iii)] There is a constant $c_3$ such that for all
$n$-dimensional
subspaces $E_1 \subset E$ one has
\[ \opo(T\iota_{E_1}) \kl c_3 \pl n^{1-\frac{1}{r}} \pl .\]
\end{enumerate}
Moreover the best constants satisfy
\[ c_1 \kl c_3 \kl c_0 \p c_2 \kl c_0^2 \pl \kla
\frac{1}{2}-\frac{1}{r} \mer^{-1} \pl c_1\pl . \]
If $E \subset B(H)$ is an operator space and $F =\min{F}$ carries
its minimal operator space structure these conditions are equivalent to
\[ \sup_{\ken} k^{1/r} \bet \la_k(TS) \rag \pl \le \pl c_4 \pl
\cgi(S)\pl  \]
for all completely $\infty$-factorable operators $S$.
\end{theorem}\hz

{\bf Proof:} \boldmath$i) \Rightarrow ii)$\unboldmath\pl By the remark
above we have for
all $S \in CB(C(K),E)$
\[ \pi_1^n(uS) \pl \le \pl c_1 \pl \noo S \rrm_{cb}\pl n^{1-1/r} \pl .
\]
By Proposition \ref{connect} this implies for all $R\in \B(F,C(K))$
\for
\sup_{\ken} k^{1/r} \bet \la_k(uSR) \rag &\le& c_0 \pl \lrx(uSR)\kl
c_0^2 \pl \kla\frac{1}{2}-\frac{1}{r}\mer^{-1} \pl c_1 \noo S \rrm_{cb}
\pl
\noo R \rrm\pl .
\mel
For the implication \boldmath$ii) \Rightarrow iii)$\unboldmath\pl
let $u: \ell_{\infty}^m \nach E_1$ be a completely bounded map and
$(y^*_k)_1^m \subset B_{Y^*}$ such that
\[ \noo Tu(e_k) \rrm \lel \langle Tu(e_k) , y^*_k \rangle \pl .\]
We define the operator $S: Y \nach \ell_{\infty}^m; S(y) \lell (\langle
y, y_k^*\rangle)_1^m$ which is of norm at most
1 and get
\for
\summ_1^m \noo Tu(e_k) \rrm &=& tr(STu) \kl 2 \pl n^{1-\frac{1}{r}}\pl
\sup_k k^{\frac{1}{r}} \bet \lambda_k(STu) \rag \kl
2 \pl  n^{1-\frac{1}{r}} c_2 \noo S \rrm \pl \noo u \rrm_{cb} \pl .
\mel
The implication \boldmath$iii) \Rightarrow ii)$\unboldmath\pl is
obvious.
Since $\lin$ is a completely complemented subspace of $M_n$ we only
have to show
the eigenvalue estimate. In fact, let $S =PR$, $R:\min(F)\nach B(H)$,
$P:B(H) \nach E$ completely bounded.
Since $F$ is considered as a subspace of $C(K)$ for some compact
Haussdorf
space $K$, there is a completely bounded extension $\hat{R}:C(K) \nach
B(H)$
of the same cb-norm by Wittstock's extension theorem, see \cite{PAU}.
If we apply $ii)$ to $S = (P\hat{R})\iota_F$, $\iota_F$ the inclusion
map we obtain the assertion.\qed

\section{1-summing operators in connection with
minimal and exact operator spaces}
\setcounter{lemma}{0}
In contrast to Banach space theory there are infinite dimensional
operator spaces such that the identity is 1-summing. This is possible
because this notion does not respect the whole operator space
structure. In fact we will see that these examples appear
in different contexts. We will start with a probabilistic approach.

\begin{lemma} \label{prob} Let $n,N \in \nz$. Then there exists a
biorthogonal
sequence $(x_j)_1^n \subset M_N$, i.e. $tr(x_j^*x_i^{ }) \lell
\delta_{ij}$ with
\[ \noo \summ_1^n e_j \otimes x_j:\ell_2^n \nach M_N\rrm_{op} \kl
\pi(1+\sqrt{2}) \pl \kla \frac{1}{\sqrt{N}} +
\frac{\sqrt{n}}{\sqrt{2}N} \mer \pl .\]
In fact a random frame for $n$-dimensional subspaces of $M_N$ satisfies
this inequality
up to a constant.
\end{lemma}

{\bf Proof:} Let $J$ be a subset of cardinality $n$ in $I \lel \{ (i,j)
\mitt i,j=1,..,N\}$.
We set $y_s\p := \p e_i \otimes e_j \in M_N$, but $x_s \p:=\p e_i
\otimes e_j$
only for $s \in J$ and $0$ else. For $(s,t)\in I\times I$ let $h_{s,t}
\p=\p \frac{1}{\sqrt{2}}(g_{st}+ig_{st}^{'})$ be
a sequence of independent, normalized, complex gaussian variables
(Clearly $(g_{st})$ and $(g^{'}_{st})$ are assumed to be independent.)
Applying Chevet's
inequality twice we obtain
\for
\lefteqn{ \ez \noo \summ_{s\in J, t \in I} h_{s,t} x_s \otimes y_t
\rrm_{op}}\\
&=& \ez \noo \summ_{s \in J,t\in I} g_{s,t} x_s \otimes
\frac{y_t}{\sqrt{2}} +
\summ_{s \in J,t\in I} g^{'}_{s,t} (ix_s) \otimes
\frac{y_t}{\sqrt{2}}\rrm\\
&\le& \kla \omega_2\{x_s,ix_s\} \pl \ez \noo \summ_{t\in I} \frac{g_t
+g_t^{'}}{\sqrt{2}} y_t \rrm_{M_N} + \frac{1}{\sqrt{2}}\p
\omega_2\{y_t,y_t\} \pl \ez \noo \summ_{s\in J} g_s x_s + g^{'}_s i x_s
\rrm_{(S_2^N)^*} \mer \\
&\le& \pl \kla 2\sqrt{N} +  \sqrt{2n}\mer \pl ,\\
\mel
where $\omega_2\{y_t,y_t\}$ corresponds to the operator norm of the
corresponding real linear operator.
Using the comparison principle between random unitary matrices in
$U_{N^2}$
and gaussian $N\times N$ matrices, see \cite{MAP}, we get
\for
\ez \noo \summ_{s \in J} x_s \otimes U(x_s) \rrm_{op}
&=& \ez \noo \summ_{s,t \in I} \langle y_t,U(x_s)\rangle x_s \otimes
y_t \rrm\\
&\le& \frac{\pi(1+\sqrt{2})}{2\sqrt{N^2}} \pl
\ez \noo \summ_{s,t} h_{s,t} x_s \otimes y_t \rrm\\
&\le& \pi(1+\sqrt{2})
\pl \kla \frac{1}{\sqrt{N}}\p + \p \frac{\sqrt{n}}{\sqrt{2}N} \mer \pl
.
\mel
For every $\eps >0$ we can find a unitary $U$ such that the norm
estimate
is satisfied up to $(1+\eps)$ by Chebychev's inequality. By passing to
a limit
we can even find a unitary $U$ satisfying the norm estimate for
$\eps=0$.
Since $U$ is a unitary in $\ell_2^{N^2}$ we use the usual
identification
between trace and scalar product to see that the elements $U(x_s)$ are
biorthogonal. An application of the concentration phenomenon
\cite{MIS} gives the assertion for random frames of $n$-dimensional
spaces
subspaces of $M_N$. \qed

The notion of random subspaces of a given $N$-dimensional Banach space
$F$
is always defined by a "natural" scalar product and the group of
unitaries of
the associated Hilbert space. A property of random subspace, means that
this
property is satisfied with "high probability" for subspaces of a fixed
dimension $n$.
In this case the probability measure is taken from the surjection $U
\mapsto span\{U(e_1),..,U(e_k)\}$
with respect to the normalized Haar measure on the group of unitaries.
Implicitly,
it is understood that the constant may depend on how close to 1 the
probability is
chosen. However, if the expected value can be estimated the
concentration phenomenon
on the group of unitaries yields reasonable estimates . For further and
more
precise information of this concept see the book of Milman/Schechtman
\cite{MIS}.
In this sense we formulate the following\hz

\begin{samepage}
\begin{cor} Let $n \kll N$ and $E$ a random subspace of $M_N$, then $E$
is 1-summing
with
\[ \opo(id_E) \kl C \pl ,\]
where $C$ depends on the probability not on the dimension.
\end{cor}
\end{samepage}

{\bf Proof:} We keep the notation from the proof above. A random
$n$-dimensional
subspace of $M_N$ is of the form $E \lell span\{U(x_s) \mitt s \in
J\}$. By lemma \ref{prob}
we can assume that with high probability the operator
\[ v \pl :=\pl \summ_{s\in J} x_s \otimes U(x_s)\]
is of norm less than  $\frac{C}{\sqrt{N}}$. The operator $vv^*$ acts as
a projection
$E$ and therefore we have the following factorization
\[ Id_E \lel (\sqrt{N}v)(\sqrt{N}v)^* \pl (\frac{1}{N} Id:M_N \nach
S_1^N) \pl \iota_E \pl ,\]
where $\iota_E$ is the canonical embedding and
$(\sqrt{N}v)(\sqrt{N}v)^*$ should
be considered as an operator from $S_1^N$ to $M_N$. As such it is of
norm at most
$C^2$. By the trivial part of the factorization theorem for 1-summing
operators \ref{fac} we get the assertion.
\qed

Paulsen, \cite{PAU}, proved that a unique operator space structure
for a given Banach space is only possible in small dimensional spaces.
This is base on the study of cb maps between minimal and maximal
operator spaces.
In this setting the author discovered lemma \ref{prob} above
in a preliminary version of this paper, noticing that this implies
an estimate for the operator integral norm for the identity
$\max(\ell_2^n) \nach \min(\ell_2^n)$.
Indeed, such a factorization has just been constructed with help of the
random spaces $E$ above.
However, the constant which can be deduced from this approach is worse
than that obtained by Paulsen/Pisier. Before we indicate
our proof of Paulsen/Pisier
result let us recall an easy lemma which is merely the definition of
the dual space, see also \cite{JP}.
For this we will use the following notation
$\noo T \rrm_n \p:=\p \noo Id_{M_n} \otimes T: M_n(E)\nach M_n(F)\rrm$
for an operator $T$
between to operator spaces $E$, $F$.

\begin{lemma} \label{ele} Let $E$, $F$ operator spaces and $T:E \nach
F$ then we have
\[\noo T \rrm_n \lel \sup \left \{ \summ_{ijkl=1}^n \langle y_{ij},
a_{ik}T(x_{kl})b_{lj} \rangle \mitt hs(a),\p hs(b) \kll 1, \pl
\noo x_{ij} \rrm_{M_n(E)},\pl \noo y_{ij} \rrm_{M_n(F^*)}\kll 1\right
\} \pl .\]
\end{lemma}

With the probabilistic approach we can prove that it
suffices to consider $n\times n$ matrices for rank $n$ operators
between minimal and maximal spaces, improving Paulsen/Pisier's result.

\begin{prop} Let $E$ be a minimal, $F$ be a maximal operator space and
$T :E \nach F$ an operator of rank at most $n$ then we have
\[ \gamma_2^*(T) \kl 170 \pl \noo Id_{M_n} \otimes T: M_n(E) \nach
M_n(F) \rrm \pl .\]
Furthermore, for every $n$ dimensional space we have
\[ \sqrt{n} \kl (\pi (1+\sqrt{2}))^2  \pl \noo Id: \min(E) \nach
\max(E) \rrm \pl ,\]
where $\min(E)$, $\max(E)$ means $E$ equipped with its minimal, maximal
operator space
structure, respectively.
\end{prop}

{\bf Proof:} First we will prove an estimate for operators $T:\ell_2^n
\nach \ell_2^n$
\[ \bet tr(T) \rag \kl 170 \noo T: \min(\ell_2^n) \nach \max(\ell_2^n)
\rrm_n \pl .\]
Indeed, we use $N\lell n$ in lemma \ref{prob} and consider the elements
\[ z_{kl} \lel \summ_i \langle x_i(e_k), e_l \rangle \otimes e_i \in
M_n(\min(\ell_2^n))\]
which are of norm at most  $\frac{\pi
(2+\frac{3}{\sqrt{2}})}{\sqrt{n}}$. In lemma \ref{ele}
we use $a=b= \frac{1}{\sqrt{n}}Id_{\ell_2^n}$ to deduce
\for
\bet tr(T) \rag &=& \bet \summ_1^n \langle T(e_i), e_j\rangle \rag \lel
 \bet \summ_{i,j} tr(x^{ }_ix_j^{*}) \langle T(e_i), e_j\rangle \rag\kl
\bet \summ_{kl} \langle z^{*}_{kl}, T(z_{kl}) \rangle \rag \\
&\le& \noo T \rrm_n \pl hs(id)^2 \pl \noo z \rrm_{M_n(\min(\ell_2^n))}
\noo z^* \rrm_{M_n(\min(\ell_2^n))}\\
&\le&
\noo T\rrm_n \pl n \pl \frac{ \pi^2 (2+\frac{3}{\sqrt{2}})^2 }{n}\kl
170 \pl \noo T \rrm_n \pl .
\mel
For an arbitrary operator $T: \min(E) \nach \max(F)$ we use trace
duality.
Indeed, let $S:F \nach E$ an operator which factors through a Hilbert
space, i.e.
$S \lel uv$, $v: F \nach H$, $u:H \nach E$. In order to estimate the
trace
we can modify $S$ by inserting the orthogonal projection on
$v(Im(T))$.
Therefore there is no loss of generality to assume $H\lel \ell_2^n$.
Hence we get
\for
\bet tr(TS) \rag &=& \bet tr(vTu) \rag\kl
170 \pl \noo vTu \rrm_n \kl
170 \pl \noo v \rrm_n \pl \noo T \rrm_n \pl \noo u \rrm_n
\kl 170 \pl \noo T \rrm_n \pl \noo v \rrm \noo u \rrm \pl .
\mel
We used that by the definition of the minimal operator space every
operator
with values in $\min(E)$ is automatically completely bounded. Taking
the infimum
over all factorizations we get the first assertion. The second one
follows
from duality by applying the estimate for the identity operator and
recalling that John's theorem $\gamma_2(Id_E) \kll \sqrt{n}$.
The better constant is obtained
by letting $N$ tend to infinity in lemma \ref{prob} and the
corresponding modification in the proof above.
\qed

As a consequence one obtains that the identity on $\max(\ell_2)$ is
indeed a
1-summing operator. More general results
hold in the context of duals
of exact operator spaces using the key inequality of \cite{JP}. We will
need some notation.
Given a Hilbert space $H$ there are at least two natural ways to
associate an operator spaces
with $H$, the column space
\[ C_H \pl:= \pl \{ x \otimes y \in B(H)\mitt x \in H \}\quad
\mbox{and the row space}\quad R_H \pl:= \pl \{ y \otimes x \in B(H)
\mitt x \in H \}\pl ,\]
where $y$ is a fixed, normalized element in $H$. It is quite easy to
check that
the corresponding norm of a matrix $(x_{ij}) \subset H$ is given by
\[\noo x_{ij} \rrm_{M_n(C_H)} \lel \noo \kla \summ_k \langle
x_{ki},x_{kj}\rangle\mer_{ij} \rrm_{M_n}^{\frac{1}{2}}
\pl \mbox{and}\pl
\noo x_{ij} \rrm_{M_n(R_H)} \lel \noo \kla \summ_k \langle
x_{jk},x_{ik}\rangle\mer_{ij} \rrm_{M_n}^{\frac{1}{2}} \pl ,\]
where we assume the scalar product $\langle \cdot,\cdot\rangle$ to be
antilinear in the first
component. It turns out that $R_H^*\lel C_H$ in the category of
operator spaces. The space
$R_H \cap C_H$ is $H$ equipped with matrix norm given by the supremum
of $R_H$ and $C_H$. The dual space $(R_H \cap C_H)^* \lell C_H +R_H$
carries
a natural operator spaces structure and was intensively studied by
Lust-Picard, Haagerup and Pisier, see \cite{LPP,HP} .
In connection with this row and column spaces it is very useful to
consider the
following notion. Let $E\subset B(K)$ an operator space and $F$ Banach
space.
An operator $T :E\nach F$ is $(2,RC)$-summing if there exists a
constant $c>0$
such that
\[ \summ_k \noo T(x_k) \rrm \kl c \pl   \max\left \{ \noo \summ_k
x_k^*x_k^{ } \rrm_{B(K)}, \noo \summ_k x_k^{ }x_k^{*} \rrm_{B(K)}
\right \}\pl . \]
The best possible constant is denoted by $\pi_{2,RC}(T)$. Let us note
that the right hand side is a weight in the sense of \cite{PSI}. We
start with a
description of $(2,RC)$ summing operators with values in a Hilbert
space,
which was suggested by C. le Merdy.

\begin{prop} \label{ext}Let $E\subset B(K)$ an operator space, $H$ a
Hilbert space
and $T:E\nach H$ a bounded linear operator. $T$ is $(2,RC)$ summing if
and only if there is a bounded extension $\hat{T}: B(K) \nach H$ of $T$
if and only if there is a completely bounded extension $\hat{T}: B(K)
\nach R_H+C_H$ of
$T$.
\end{prop}

{\bf Proof:} Let us observe that by the non-commutative Grothendieck
inequality
every bounded $S:B(K) \nach H$ is $(2,RC)$ summing, see e.g.
\cite{PIL}. Therefore,
we are left to prove the existence of a cb extension for
$(2,RC)$-summing operators
$T:E \nach H$. Using  a variant of Pietch's factorization theorem,(for
more
precise information see \cite{PSI},) there are states $\phi$, $\psi$ on
$B(K)$ and $0\le \theta \le 1$ such that
\[ \noo v(x) \rrm \kl \pi_{2,RC}(T) \pl \kla \theta \phi(xx^*) +
(1-\theta) \psi(x^*x) \mer^{\frac{1}{2}} \pl .\]
We define the sesquilinearforms $\langle x,y \rangle_{\phi} \p:=\p
\phi(yx^*)$ and
$\langle x,y \rangle_{\psi} \pl :=\pl \psi(x^*y)$.
Furthermore, we denote by $C_{\phi}$, $R_{\psi}$ the column, row
Hilbert space which is induced
by the corresponding scalar product. It is easy to check that the
identities
$I_{\phi}: B(K) \nach C_{\phi}$, $I_{\psi} :B(K) \nach R_{\psi}$ are in
fact
completely bounded of norm $1$. We denote by
\[ M\pl:=\pl cl\{ (\sqrt{\theta}x,\sqrt{1-\theta}x) \mitt x \in E \}\pl
\subset \pl C_{\phi} \oplus_2 R_{\psi}\]
the closure of the image of $J \p:=\p\sqrt{\theta}I_{\phi} \oplus
\sqrt{1-\theta}I_{\psi}$ restricted
to $E$. $P_M$ denotes the orthogonal projection of $C_{\phi}\oplus_2
R_{\psi}$ onto $M$.
Then we get an extension  $\hat{T} \lel \tilde{v}J$ of $T$, where
$\tilde{v}$ acts as $v$ but considered as an operator on $M$.
By the first inequality $\tilde{v}$ is of norm at most $\pi_{2,RC}(T)$
and by definition
of $R_H+C_H$ we get $\noo \tilde{v}:R_M+C_M \nach R_H+C_H\rrm_{cb}
\kll  \pi_{2,RC}(T)$.
By duality it is easy to see that $P_M: C_{\phi}\oplus_1 R_{\psi} \nach
R_M+C_M$ is completely
bounded of norm 1. On the other hand the cb norm of
\[  \sqrt{\theta}Id_{C_{\phi}} \oplus \sqrt{1-\theta}Id_{R_{\psi}}:
C_{\phi}\oplus_{\infty} R_{\psi} \nach C_{\phi}\oplus_1 R_{\psi} \]
is at most $\sqrt{2}$. \qed

Now we will give a description of completely bounded operators between
the class of exact operator spaces and maximal operator spaces.
Pisier's notion of exact operator spaces, \cite{PSE}, is motivated
by Kirchberg's work.
One possible definition says that an operator space
is exact if its finite dimensional subspaces are uniformly cb
isomorphic
to subspaces of the spaces of compact operators.

\begin{samepage}
\begin{prop} \label{cb}
Let $E \subset B(K)$ be either an exact operator space and $F$ a
maximal operator space, i.e, a quotient
of $\ell_1(I)$ for some index set $I$,
or $E$ a $C^*$ algebra and $F \lell \ell_1(I)$.
For an operator $T:E \nach F$ the following are equivalent.
\begin{enumerate}
\item[i)] $T$ is completely bounded.
\item[ii)] There is a  Hilbert space $H$ and operators $v:E \nach H$,
$u: H \nach F$ such that
$v$ is $(2,RC)$ summing and $u^*$ is absolutely $2$ summing.
\item[iii)] There is a completely bounded extension $\hat{T}: B(K)
\nach \ell_1(I)$
of $T$.
\item[iv)] $T$  factors completely through $R_H + C_H$ for some Hilbert
space $H$.
\end{enumerate}
Moreover, the corresponding constants are equivalent.
\end{prop}
\end{samepage}
{\bf Proof:}  The implication $i) \Rightarrow ii)$ is either the
non-commutative
Grothendieck inequality, see \cite{PIL}, or the key inequality  in
\cite{JP}.
The implications $ii) \Rightarrow iii), iv)$ are direct consequences of
proposition \ref{ext} and
the extension properties of absolutely 2 summing operators. We only
have to note
that an absolutely $2$ summing operator $u^*:\ell_{\infty}(I)\nach R_H
\cap C_H$ is
completely bounded. The rest is trivial. \qed

For the proof of theorem 3 we will need some more notation. Let $1 < p
< \infty$, $E$ be an operator space
and $F$ a Banach space. An operator $T:E \nach F$ belongs to
$\Gamma_{p,RC}$ if
\[ \grc(T) \pl :=\pl \sup\{ \si_{p,\infty}(vTu) \mitt v\in
\Pi_2(F,\ell_2), \pl u \in CB(R+C,E)\pl
\pi_2(v), \pl \noo u\rrm_{cb}\kl 1\} \pl < \pl \infty\pl .\]
For $p=1$ we will use $\Gamma_{1,RC}$, $\gamma_{1,RC}$ for the
corresponding expressions
with $\si_{p,\infty}$ replaced by $\si_1$.
This notion is modeled close to the notion of Hilbert space factoring
operators
and forms a 'graduation' of $\Gamma_{1,RC}$ in the cases $p>1$.
This has already been proved to be useful for
eigenvalue estimates.

\begin{theorem} \label{exact} Let $1< p< 2$, $G$ an exact operator
space, $E\subset G^*$ and $F$ a minimal operator space.
For an operator $T:E \nach F$ the following are equivalent.
\begin{enumerate}
\item[i)] There exists a constant $c_1>0$ such that
\[ \summ_1^n \noo Tx_k  \rrm \kl c_1 \pl n^{1-\frac{1}{p}} \pl \noo
\summ_1^n e_i \otimes x_i \rrm_{\ell_1^n \otimes_{min} E} \pl. \]
\item[ii)] $T$ is in $\Grc$
\item[iii)] There is a constant $c_3$ such that for all completely
bounded operator $S:F \nach E$ one has
\[ \sup_k k^{\frac{1}{p}} \pl \bet \lambda_k(TS) \rag
\kl c_3 \pl \noo S:{\rm min}(F) \nach E \rrm_{cb} \pl .\]
\end{enumerate}
In the limit cases $p=1$ the same remains true if we replace the
$\ell_{p,\infty}$ norm of the eigenvalues by the $\ell_1$ norm.
Furthermore, every completely bounded $S: \min{F} \nach E$ is
absolutely 2-summing
and hence the eigenvalues of a composition $TS$ are in $\ell_2$.
\end{theorem}

{\bf Proof:} The implication $iii) \Rightarrow i)$ follows along the
same line as $iii)\Rightarrow i)$ in \ref{opeigen}.
For the implication \boldmath $ii) \Rightarrow iii)$ \unboldmath let
$S: \min(F) \nach E$ completely bounded
and consider $\hat{S} \p:=\p\iota_{E}S: \min(F) \nach G^*$. By
proposition \ref{cb} we can assume that
$\hat{S} \lell uv$, where $v:F \nach H$ is absolutely 2-summing and
$u:R_H \cap C_H \nach G^*$
is completely bounded. Using an orthogonal projection $P$ on
$u^{-1}(E)$ together
with the homogeneity of the spaces $R_H$ and $C_H$, see \cite{BPT}, we
can assume that
$u(H) \subset E$ and therefore $S\lell uv$. Using the well-known
eigenvalue estimate of the class
$\es_{p,\infty}$ and the principle of related operators, \cite{PIE},
we get
\for
\sup_k k^{\frac{1}{p}} \pl \bet \lambda_k(TS) \rag &=&
\sup_k k^{\frac{1}{p}} \pl \bet \lambda_k(vTu) \rag \kl
 c_0 \pl\sup_k k^{\frac{1}{p}} \pl a_k(vTu) \\
&\le& c_0 \pl \grc(T) \pl \noo u \rrm_{cb} \pl \pi_2(v) \kl
c_0 \pl b_0 \pl \noo S \rrm_{cb} \pl \grc(T) \pl ,
\mel
where $b_0\kll 4 \sqrt{2}$ is the constant from proposition \ref{cb}.
In order to prove \boldmath $i) \Rightarrow ii)$ \unboldmath we will
use
the notion of Grothendieck numbers for an operator $R: X \nach Y$
introduced by S. Geiss.
\[ \Gamma_n(R) \pl :=\pl \sup \left\{ \bet {\rm det}(\langle R(x_i),
y_j \rangle)_{ij} \rag^{\frac{1}{n}} \mitt
(x_i)_1^n \subset B_X, \p (y_j)_1^n\subset B_{Y^*} \right \} \pl. \]
Using an inequality of  \cite{DJ1} we have to show that
\[\sup_n n^{\frac{1}{p}-\frac{1}{2}}\pl \Gamma_n(Tu) \kl c_2 \pl \noo u
\rrm_{cb} \]
for all operator $u : R+C \nach E$. By the definition of the
Grothendieck numbers
we have to consider elements $(y_k^*)_1^n \subset B_{F^*}$ and $v:=
\summ_1^n y_k^{*}\otimes e_k:F \nach \ell_{\infty}^n$
which is of norm $1$. If $\iota_{2,\infty}^n :\ell_{\infty}^n \nach
\ell_2^n$
denotes the canonical inclusion map we have to show
\[ \Gamma_n(\iota_{2,\infty}^n vTu) \kl c_2 \pl
n^{\frac{1}{2}-\frac{1}{p}} \pl \noo u \rrm_{cb} \pl .\]
Now let  $w: \ell_2^n \nach H$ such that
\[ \summ_1^n a_j(\ioc vTu) \lel tr(\ioc vTuw) \pl .\]
Using basic properties of Grothendieck numbers, see \cite{GEI} and the
geometric/arithmetic mean
inequality we get for $S:= uw\ioc : \ell_{\infty}^n \nach E$
\for
\Gamma_n(\iota_{2,\infty}^n vTu) &\le&
\kla \prod_1^n a_j(\ioc vTu) \mer^{\frac{1}{n}}
\kl \frac{1}{n} \pl \summ_1^n a_j(\ioc vTuw) \lel
\frac{1}{n} \pl \bet tr(vTS) \rag \\
&\le& \frac{1}{n} \pl \summ_1^n \noo TS(e_k) \rrm \pl \sup_{i} \noo
y_k^* \rrm \kl
c_1 \pl n^{-\frac{1}{p}} \pl \noo S \rrm_{cb} \\
&\le& c_1 \pl n^{-\frac{1}{p}} \pl \pi_2(\ioc) \pl \noo wu :R_n\cap C_n
\nach E\rrm_{cb}\kl
c_1 \pl n^{\frac{1}{2}-\frac{1}{p}}\pl \noo u \rrm_{cb}\pl ,
\mel
where we have used the homogeneity of the space $R_H\cap C_H$ and
remark \ref{rcoh} to estimate
the cb norm of $\ioc: \ell_{\infty}^n \nach R_n\cap C_n$.
In the case $p=1$ we have to estimate $\si_1(vTu)$ for a 1-summing $T$,
absolutely
2 summing v and completely bounded $u:R_n\cap C_n \nach E$
By Pietsch's factorization theorem, see \cite{PIE}, there is
a factorization of $v \lell SR$, $S:\ell_{\infty}^M \nach \ell_2$, with
absolutely 2 summing $S$. Since $uwS$ is completely bounded
for all bounded $w$ we can use the definition  to see that $TuwS$ is
integral
in the Banach space sense and hence the trace of $wvTu \lell SRTuw$ can
be estimated by
the 1-summing norm. This gives the right estimate for trace class norm,
and hence the eigenvalues of $ST$, provided
$w$ is chosen by polar decomposition as above. \qed

\begin{rem} \label{fact}{\rm  A variant of Kwapien theorem for Hilbert
space factorizing operators
shows that an operator $T: G^* \nach \min(F^{**})$
factors completely through $R_H\cap C_H$ if and only if $|tr(TS)|\kl C$
for all operators $S:F^{**}\nach G^*$ which admit a factorization
$S=vu$, $\pi_2(v)\le 1$ and $\pi_{2,RC}(u^*)\le 1$. Indeed, this
duality
concept was studied in the more general framework of $\gamma$-norms by
Pisier \cite{PSI}. We want to indicate the connection to 1-summing
operators in this
context. Given a 1-summing operator $T:E\subset G^* \nach F$ we observe
that
$T$ corresponds by trace duality to a linear functional
on $F^*\otimes_{\alpha} E$ where $\alpha(S) \p:=\p \inf\{\pi_2(v) \p
\noo u :R \cap C \nach E \rrm_{cb}\}$
and the infimum is taken over all factorizations $S=vu$. Since
$F^*\otimes_{\alpha}E$ embeds isometrically into $F^* \otimes_{\alpha}
G^*$
an application of Hahn-Banach yields a norm preserving functional on
the whole tensor product,
i.e. an extension $\hat{T}:G^* \nach F^{**}$ of $T$, which is also
1-summing by theorem \ref{exact}.
As a consequence of the key inequality in \cite{JP} and
proposition \ref{ext} we deduce that for all $u:R_H\cap C_H \nach G^*$
the cb-norm
is equivalent to $\pi_{2,RC}(u^*)$. Therefore, we can apply the
modification of Kwapien's argument, see also  \cite{PSI}, to obtain a
completely bounded
factorization of $\hat{T}: G^* \nach \min(F^{**})$
through $R_H\cap C_H$ for some Hilbert space $H$. Clearly, if $\hat{T}$
admits such a factorization it must be 1-summing and all these
properties coincide
due to the fact that $G$ is exact.}
\end{rem}

\begin{cor} Let $G$ be an exact operator space, $q:B(H)^* \nach G^*$
the quotient map
and $E \subset G^*$.
The following conditions are equivalent
\begin{enumerate}
\item The Banach space $E$ is of cotype 2 and every bounded operator
$u:c_0 \nach E$ is completely bounded.
\item The Banach space $E$ is of cotype 2 and every operator $v: E
\nach R\cap C$
which admits a completely bounded extension $\hat{v}:G^* \nach R\cap C$
is absolutely 1-summing.
\item There exists a constant $c>0$, such that for every sequence
$(x_k)_1^n\subset E$ there is a sequence $(\tilde{x}_k)_1^n \subset
B(H)^*$
such that $q(\tilde{x}_k)=x_k$ and
\[ \ez \noo \summ_1^n \tilde{x}_k \eps_k \rrm_{B(H)^*} \kl c \pl \ez
\noo \summ_1^n x_k \eps_k \rrm_{E} \pl .\]
\end{enumerate}
In particular, a maximal operator space satisfies one of the conditions
above
if and only if it is of operator cotype 2 if and only if it is a G.T.
space of
cotype 2, see \cite{PIL}.
\end{cor}

{\bf Proof:} Let $X$ be a Banach space we define $Rad(X) \subset
L_2(\dz,X)$
to be the span of $\{ \eps_i \otimes x_i \}$, where
$\dz =\{-1,1\}^{\nz}$ is the group of signs with its Haar measure $\mu$
and $\eps_i$ the i-th coordinate.
For a sequence
$(x_i)_i$ the norm in $Rad(X)$ is given by
\[ \noo (x_i)_i \rrm \pl :=\pl \kla \int_{\dz} \noo \summ_i \eps_i x_i
\rrm_{X}^2 d \mu \mer^{\frac{1}{2}} \pl .\]
It was shown by Pisier and Lust-Picard \cite{LPP} that
$Rad(B(H)^*)$ and $(R+C)(B(H)^*)$ have equivalent norms. Since the map
$Id_{R+C}\otimes q$ is a complete quotient map, condition $iii)$ is
equivalent to
\[ \noo Id\otimes \iota_E: Rad(E) \nach (R+C)(G^*) \rrm \pl < \pl
\infty \pl ,\]
where $\iota_E: E \nach G^*$ is the inclusion map. We deduce
from theorem \ref{exact} and remark \ref{fact} that  condition $i)$ and
$ii)$ are
equivalent by trace duality. Moreover, all conditions imply that $E$ is
of cotype 2,
since $B(H)^*$ is of cotype 2, \cite{TOJ}. Now let $v : = \summ_i x_i^*
\otimes e_i$
be an operator from $E$ to $R\cap C$. We deduce from \cite[5.16]{PIL}
\[ \frac{1}{C_1(E)}\pl \pi_1(v) \kl \pi_2(v) \kl \noo (x^*_i)^{ }_i
\rrm_{(Rad(E))^*}
\kl C_2(E) \pl \pi_2(v) \kl C_2(E) \pl \pi_1(v) \pl ,\]
where $C_2(E)$ is the cotype 2 constant of $E$ and $C_1(G)$ only
depends
on $C_2(E)$. Finally we note that $CB(G^*,R\cap C) \lel (R\cap
C)(G^{**})$.
But this means that the set of operators admitting a cb extension can
be identified
with the dual space of $(R+C)^{inj}(E) := (Id \otimes
\iota_E)^{-1}(R+C)(G^*)$.
Therefore condition $ii)$ is equivalent to
\[ \noo Id \otimes Id_G:(R+C)^{inj}(E) \nach (Rad(E))^* \rrm <
\infty\pl.\]
Duality implies the assertion. In the situation of maximal operator
spaces
we deduce from remark \ref{ocot} that a maximal operator space
$X=\ell_1(I)/S$
with operator cotype 2 satisfies condition $i)$ whereas $iii)$ implies
operator
cotype 2 since $\ell_1(I)$ has operator cotype 2. (This seems not to be
the case for
$\es_1(H)$.)\qed

In the last part we will study the operator spaces associated to
the Clifford algebra. Recalling that the generators of the Clifford
algebra
have already been useful to find an example of a $(2,oh)$-summing
space, see \cite{PLT},
it is probably not surprising that this space is also 1-summing. More
precisely,
let $(u_i)_{i\in N} \subset {\otimes \atop \nen}M_2$ be the generators
of the
Clifford algebra, i.e.
\for
 u_i \lel u_i^* &\pl \mbox{and}\pl& u_i^2 \lel Id \quad \quad
 \mbox{for} \pl i\in \nz\pl ,\\
 u_iu_j\p +\p u_ju_i &=& 0 \quad \quad \mbox{if}\pla i \neq j \pl.
\mel
By $CL$ we denote of the span of these generators. The next proposition
collects some
facts about this space. ($OH$ is the operator Hilbert space introduced
and studied by Pisier \cite{PLT}).

\begin{prop}\label{Sn}
\begin{enumerate}
\item $CL$ is $\sqrt{2}$ isomorphic to a Hilbert space.
\item The identity $id_{CL}$ is 1-summing with $\opo(T) \kll 2$ and for
every operator $T:CL\nach CL$ we have
\[ \summ_k \bet \la_k(T) \rag \kl c_0  \pl \cgi(T) \pl . \]
\item Let $G\in \{ OH,C,R,C+R,R\cap C\}$ and $u:G \nach CL$ then one
has
\[ \noo u \rrm_{cb} \sim_c \pi_2(u) \pl .\]
\end{enumerate}
\end{prop}

{\bf Proof:} By approximation it is sufficient to consider the finite
dimensional case. Therefore we fix
$\nen$ and $u_1,..,u_n \in \otimes_{k=1}^n M_2 \cong M_{2^n}$. For an
element
$x \lell\summ_j \alpha_j u_j$  we have
\for
x^*x+xx^* &=& \summ_{kj} \overline{\alpha_k} \p \alpha_j\p  u_k u_j\pl
+\pl
\summ_{jk} \alpha_j\p \overline{\alpha_k} \p u_j u_k\\
&=& 2 \pl\summ_1^n \bet \alpha_k \rag^2 u_k^2 \pl  +\pl \summ_{k<j}
\overline{\alpha_k}\p \alpha_j\p  (u_ku_j+u_ju_k)
+\pl+\pl \summ_{k>j} \overline{\alpha_k}\p \alpha_j\p
(u_ku_j+u_ju_k)\\
&=& 2 \pl \noo \alpha \rrm_2^2 \pl Id \pl .
\mel
In particular, we get
\[\noo \alpha \rrm_2 \lel \noo \frac{x^*x+xx^*}{2}
\rrm_{M_{2^n}}^{\frac{1}{2}}
\kl \pl \noo x \rrm_{M_{2^n}} \kl \sqrt{2} \noo \frac{x^*x+xx^*}{2}
\rrm^{\frac{1}{2}}
\lel \sqrt{2} \noo \alpha \rrm_2 \pl. \]
This is the first assertion.  In order to estimate the 1-summing norm
we define $\hat{x} \lell \kla {x \atop x^*} {0 \atop 0} \mer$
in $\es_1^{2^{n+1}}$. With the triangle inequality in $\es_1^{2^{n+1}}$
we get
\[2^n \pl \noo \alpha \rrm_2 \lel \noo \kla \frac{x^*x+xx^*}{2}
\mer^{\frac{1}{2}} \rrm_1
\lel \frac{1}{\sqrt{2}} \pl \noo \hat{x} \rrm_{\es_1^{2^{n+1}}} \kl
\sqrt{2} \noo x \rrm_{\es_1^{2^n}} \pl .\]
Combining these estimates we have found a factorization of the identity
on $CL^n$ through
the restriction of $2^{-n} Id:M_{2^n} \nach \es_1^{2^n}$ on $E$.
By proposition \ref{fac} the 1-summing norm of identity on $CL^n$ is
at most $2$.  As a consequence every operator $T: CL \nach CL$ which
factors
completely through a $C(K)$ space is integral and since $CL$ is
isomorphic to a Hilbert space
the eigenvalues are absolutely summing. To prove $3.$ let $u:R\cap
C\nach CL$. In order
to show that this operator is absolutely 2-summing we use trace
duality.
For this let $v:CL \nach R\cap C$ which is absolutely 2-summing. By
Pietsch factorization
theorem $v$ factors through a 2-summing operator $S:C(K)\nach R\cap C$,
which is completely bounded,
see \ref{rcoh}.
Since $CL$ is 1-summing the composition $Su$ is integral we get the
right estimate for the
trace. Vice versa, we consider an absolutely 2 summing operator $u:R+C
\nach CL$.
All the underlying Banach spaces are isomorphic to Hilbert spaces and
therefore
$u$ admits a factorization $u\lell wv$, $v^*$ absolutely 2-summing and
$w: \ell_1 \nach CL$.
This operator $w$ is automatically completely bounded, whereas $w$ is
completely bounded
in view of \ref{rcoh} and duality. \qed

Now we will construct operator spaces $E_r$, $E_r^n$ which are
isomorphic to $\ell_2$, $\ell_2^n$, respectively,
but the 1 summing norm has a certain growth rate.
For $1<r<2$ we define a matrix structure on $\ell_2$, ($\ell_2^n$) as
follows
\[ \noo (x_{ij}) \rrm_r \pl := \pl \sup_{\ken} k^{\frac{1}{r}-1}\p
\sup \left\{ \noo (P_H\p x_{ij}) \rrm_{CL} \mitt H \subset \ell_2,\p
(H\subset \ell_2^n)\p dimH
\p\le\p k  \right\}\pl,\]
where we identify $CL$ and $\ell_2$ via the isomorphism from
proposition
\ref{Sn} and $P_H \p:\p \ell_2 \nach H$ denotes the orthogonal
projection on $H$.
The next proposition states the properties of this operator spaces.

\begin{prop} Let $1<r<2<p<\infty$ with
$\frac{1}{r}=\frac{1}{2}+\frac{1}{p}$.
\begin{enumerate}
\item[i)] $E_r$ is an operator space which is $2$ isomorphic to
$\ell_2$.
\item[ii)] For all $\nen$ one has $\opo^n(id_{E_r}) \pl \sim_2 \pl
n^{\frac{1}{r}}\pl. $
\item[iii)] For all completely $\infty$-factorable operators $T \in
\Cgi(E_r,E_r)$
one has
\[ \sup_{\nen} n^{\frac{1}{r}} \p\bet\la_n(T)\rag \pl\le \pl c_r\pl
\cgi(T:\min(E_r) \nach E_r)\pl.\]
\item[iv)] For the completely bounded operators with values
in $E_r$ and defined on $\ell_{\infty}$ or $G \in \{R,C,R+C,R\cap C,
OH\}$
one has
\[ CB(\ell_{\infty},E_r) \pl=\pl \Lpa(\ell_{\infty},E_r)
\pll \mbox{and} \pll CB(G,E_r) \pl=\pl \Lpa(OH,E_r)\pl.\]
\end{enumerate}
A similar statement holds uniformly in $n$ for the spaces $E_r^n$.
\end{prop}
\pagebreak

{\bf Proof:}  i) is clear by definition and proposition \ref{Sn}.
ii) and iii) follows from iv) and standard estimates  of $\lrx(id:
\ell_{\infty}^n\nach \ell_2^n) \sim_{c_r} n^{\frac{1}{r}}$.
For $iv)$ we note that by definition and the fact that $CL$ is 1
summing we have
\[ \noo T: \ell_{\infty} \nach CL \rrm_{cb} \sim_{2} \sup_{\ken,dim H
\le k} k^{\frac{1}{r}-1} \pl \iota(P_HT)\pl .\]
For an operator $u: \ell_2 \nach \ell_{\infty}$ we deduce by Schmidt
decomposition
\[ \sup_k k^{\frac{1}{r}} a_k(Tu) \kl \sup_{\ken, dim H\le k}
k^{\frac{1}{r}-1} \si_1(P_HTu)
\kl 2 \noo T\rrm_{cb} \pl \noo u\rrm \pl .\]
For the converse implication we use an interpolation argument. Indeed,
by
standard relations between different s-numbers, \cite{PIE}, one has
\[\Lrx \pl \subset \pl \Lpa \pl \subset \pl (\L_{2,1}^{(a)},
\L_{\infty}^{(a)})_{\theta,\infty} \pl ,\]
with $\frac{1}{p}=\frac{1-\theta}{2}+\frac{\theta}{\infty}$ and
$\lze(T)$
is the norm of the approximation numbers in the Lorentz spaces
$\ell_{2,1}$.
By definition of the $K_t$ functional for $t=\sqrt{k}$ we can find a
decomposition
$T \lell T_1 +T_2$ such that $\lze(T_1) + \sqrt{k} \p\noo T_2 \rrm \kll
c_p \pl k^{\frac{1}{2}-\frac{1}{p}}\p\lrx(T)$.
An application of "little Grothendiek's inequality", \cite{PIL}, gives
$\iota(T_1) \kl c_1 \p \lze(T_1)$. Hence we get
for every $k$ dimensional subspace $H$
\for
 \iota(P_HT) &\le& \iota(P_HT_1) + \iota(P_HT_2) \kl
 \iota(T) +\sqrt{k}\p \pi_2(T_2)\\
&\le&
(c_1 + \frac{2}{\sqrt{\pi}}) \pl \kla \lze (T_1) +\sqrt{k} \p \noo P_H
T_2 \rrm \mer
\kl c_p \p (c_1 + \frac{2}{\sqrt{\pi}}) \pl k^{1-\frac{1}{r}}\pl
\lrx(T) \pl .
\mel
The second formula is proved along the same line, although
Grothendiek's inequality is not
used in this argument. The key point here is the following formula
which we deduce from
proposition \ref{Sn}
\for
 \noo T\p: \p G\nach E_r \rrm_{cb}&=&
\sup_{\ken} k^{\frac{1}{r}-1} \pl \sup_{H,\p dim(H)\le k} \noo P_HT \p:
\p G
\nach CL\rrm_{cb}\\
&\sim_c& \sup_{\ken} k^{\frac{1}{r}-1} \pl \sup_{H,\p dim(H)\le k}
\pi_2 (P_HT \p: \p G
\nach CL)\pl .\\[-1.3cm]
\mel
\qed

\begin{rem} {\rm An easy modification of the spaces above allows us to
construct
an operator space $E_1$ such that the identity is not
1-summing, but
\[ \sup_{\nen} n \pl\bet \la_n(T) \rag \pl \le \pl c_0 \pl \cgi(T) \]
for $T \in \Cgi(E_1,E_1)$. In fact, we define the matrix norm on $E_1$
by
\[ \noo (x_{ij}) \rrm \pl := \pl \sup_{\ken} \p
\sup \left\{ \noo (P_H\p x_{ij}) \rrm_{CL} \mitt H \subset \ell_2, \p
dim(H)
\p\le\p k \pl \mbox{and} \pl H \subset H_k
 \right\}\pl,\]
where $H_k\p :=\p spann\{e_j\p |\p k\p\le\p j\}$. Using similar
arguments as above we
can prove
\[ CB(\ell_{\infty},E_1) \pl\subset\pl {\cal L}_{1,
\infty}^{(x)}(\ell_{\infty},E_1) \pl, \]
and the diagonal operator $\ds \in \com(\ell_{\infty},E_1)$ defined by
$\sigma_k \p=\p \frac{1}{k}$ is completely bounded, but not 1-summing.}
\end{rem}


\begin{exam} {\rm At this point we want to give a review
of infinite dimensional operator spaces such that $\opo^n(id_E) \kl
n^{1-\frac{1}{r}}$
for some $1<r<2$.
By theorem \ref{exact} it is easy to see that this holds for $\max(X)$
if and only
if $X$ is a so called weak $r$ Hilbertian Banach spaces, see
\cite{PI3}, \cite{GEI} and \cite{DJ1}.
Standard examples are obtained by interpolation $X\lell(H,Y)_{\theta}$,
$\frac{1}{r} \lell \frac{1-\theta}{1} + \frac{\theta}{2}$,
where $H$ is a Hilbert space and $Y$ an arbitrary Banach space.
Therefore, $\max(\ell_r)$ and $\max(\ell_{r'})$ are typical examples,
but
also $\max(\es_r)$ and $\max(\es_{r'})$. Moreover, we see that
$\opo^n(Id_{\max(E)}) \kll n^{1-\frac{1}{r}}$
if and only if the same holds for $\max(E^*)$. In the limit case $r=1$
the identity
of a maximal operator space is 1-summing if and only if the associated
Banach space is isomorphic
to a Hilbert space, whereas a subspace of $\max(E)$ is 1-summing if and
only
if it is a complemented Hilbert space in $E$.

It is easy to see that the operator space $CL$ spanned by the
generators of
the Clifford algebra is an exact operator space. Moreover, the
exactness constant \cite{PSE}
is uniformly in $n$ bounded for the spaces $E_r^n$.
Using theorem \ref{exact} and the last proposition it is quite
standard to deduce that the operator space dual $CL^*$ is not
1-summing,
but $\opo^k(id_{(E_r^n)^*}) \sim_{c_r} k^{\frac{1}{r}-\frac{1}{2}}$
for $k \le n$. }
\end{exam}
\newcommand{\ww}{\vspace{-0.2cm}}

\hz
1991 Mathematics Subject Classification: 47C15, 47B06.
\begin{quote}
\newcommand{\www}{\vspace{-0.15cm}}
Marius Junge\www

Mathematisches Seminar der Universit\H{a}t Kiel\www

Ludewig-Meyn-Str. 4 \www

24098 Kiel \www

Germany\www

Email: nms06@rz.uni-kiel.d400.de\www
\end{quote}
\end{document}